\documentclass{amsart}
\usepackage[all]{xy}
\usepackage{amsfonts,amscd,amssymb,amsmath,mathrsfs}
\theoremstyle{plain}
\theoremstyle{plain}
\newtheorem{thm}{Theorem}[section]
\newtheorem{lem}[thm]{Lemma}
\newtheorem{prop}[thm]{Proposition}
\newtheorem{cor}[thm]{Corollary}
\theoremstyle{definition}
\newtheorem{defn}{Definition}[section]

\theoremstyle{remark}
\newtheorem{rmk}{Remark}[section]

\numberwithin{equation}{section}

\begin{document}

\newcommand{\Soc}{\operatorname{Soc}}
\newcommand{\modP}{\mod_{P}\Lambda}
\newcommand{\modl}{\mod \Lambda}
\newcommand{\module}{\operatorname{mod}}
\newcommand{\p}{\operatorname{p}}
\newcommand{\inj}{\operatorname{i}}
\newcommand{\Module}{\operatorname{Mod}}
\newcommand{\Cok}{\operatorname{Coker}}
\newcommand{\Hom}{\operatorname{Hom}}
\newcommand{\h}{\operatorname{h}}
\newcommand{\e}{\operatorname{e}}
\newcommand{\res}{\operatorname{res}}
\newcommand{\Tr}{\operatorname{Tr}}
\newcommand{\TrD}{\operatorname{TrD}}
\newcommand{\rad}{\operatorname {{\bold r}}}
\newcommand{\La}{\operatorname{\Lambda}}
\newcommand\End{\operatorname{End}}
\newcommand\Ext{\operatorname{Ext^1_\Lambda}}
\newcommand\Ex{\operatorname{Ext}}
\newcommand\ann{\operatorname{ann}}
\newcommand\coend{\operatorname{coend}}
\newcommand\Img{\operatorname{Im}}
\newcommand\D{\operatorname{D}}
\newcommand\DTr{\operatorname{DTr}}
\newcommand\Ker{\operatorname{Ker}}
\newcommand\Coker{\operatorname{Coker}}
\newcommand{\tr}{\operatorname {t}}
\newcommand{\ct}{\operatorname{ct}}
\newcommand{\rej}{\operatorname{rej}}
\newcommand{\gd}{\operatorname{{gl.dim}}}
\newcommand{\radic}{\operatorname{rad}}
\newcommand{\add}{\operatorname{add}}
\newcommand{\Ind}{\operatorname{Ind}}
\newcommand{\Sub}{\operatorname{Sub}}
\newcommand{\Fac}{\operatorname{Fac}}
\newcommand{\lap}{\operatorname{l}}
\newcommand{\rap}{\operatorname{r}}
\newcommand{\locn}{\operatorname{lnRep}}
\newcommand{\rep}{\operatorname{Rep}}

\newcommand{\spec}{(\Gamma,
\Lambda,\textrm{\textbf{d}},\mathfrak B)}
\newcommand{\newspec}{(\Gamma,
\sigma_x\Lambda)}
\newcommand{\gl}{(\Gamma,\Lambda)}
\newcommand{\poset}{(\Gamma_{0}¥,\Lambda)}
\newcommand{\Stilde}{\tilde S}
\newcommand{\Ttilde}{\tilde T}
\newcommand{\ad}{(+){\rm -admissible}}
\newcommand{\pa}{{\tilde{\mathcal P}}(A)}
\newcommand{\modA}{{\rm mod}A}
\newcommand{\hull}{H_{\Lambda}}
\newcommand{\modr}{\mathrm{mod}\,R_{\gl}}

\bibliographystyle{plain}
\title[Sequences of reflection functors]{Sequences of reflection functors and the preprojective component of a valued quiver}
\author{Mark Kleiner}
\address{Department of Mathematics, Syracuse University, Syracuse, New
York 13244-1150}
\email{mkleiner@syr.edu}
\author{Helene R. Tyler}
\address{Department of Mathematics and Computer Science, Manhattan
College, Riverdale, New York 10471}
\email{helene.tyler@manhattan.edu}
\dedicatory{Dedicated to the
memory of L. Gaunce Lewis, Jr.}
\thanks{The first-named author is supported by the NSA grant H98230-06-1-0043.  The paper was written when the second-named
author visited Syracuse University in June of 2004 and in July of
2006 with a partial support of a Summer Research Grant
from Manhattan College. She expresses her sincere gratitude to the
college for their generosity and to the members of the
Syracuse University Mathematics Department for their warm
hospitality.}


\begin{abstract} This paper concerns preprojective  representations of a finite connected valued quiver without oriented cycles.  For each such representation, an explicit formula in terms of the geometry of the quiver gives a unique, up to a certain
equivalence, shortest (+)-admissible sequence such that the corresponding composition of
reflection functors annihilates the representation. The set of equivalence classes of the above sequences is a partially ordered set that contains a great deal of information about the preprojective
component of the Auslander-Reiten quiver. The results apply to the study of reduced words in the Weyl group associated to an indecomposable symmetrizable generalized Cartan matrix.
\end{abstract}
\maketitle

\section*{Introduction} \label{S:intro}

The motivation for this work comes from two sources.  The first is the paper ~\cite{kt}, which assigns a canonical (+)-admissible sequence to each indecomposable preprojective module  over the path algebra of a finite connected quiver without oriented cycles and then uses the combinatorial structure of the set $\mathfrak S$ of (+)-admissible sequences, and reflection functors instead of the Coxeter functor (Auslander-Reiten translation), to give an explicit description of the preprojective component of the Auslander-Reiten quiver \cite{ars}.  In this connection a question is whether similar results hold in a more general setting of representations of valued quivers studied in ~\cite{dr}.  The question is especially relevant in view of  ~\cite{kp}, which is our second source of motivation.  Using combinatorics  of the set $\mathfrak S$, the latter paper  relates properties of reduced words in the Weyl group $\mathcal {W}(A)$ associated to an indecomposable symmetric generalized $n\times n$ Cartan matrix  $A$ ~\cite{kac1990}   to properties of  preprojective  modules over the path algebra of a quiver without oriented cycles whose underlying graph is the graph associated to $A$.  Let $\sigma_1,\dots,\sigma_n$ be the simple reflections, and let  $c$ be any Coxeter element, i.e., $c=\sigma_{x_n}\dots\sigma_{x_1}$  where $x_1,\dots,x_n$ is any permutation of the  numbers $1,\dots,n$.  The authors of  ~\cite{kp} proved that $\mathcal {W}(A)$ is infinite if and only if the powers of $c$ are reduced words in the $\sigma_h$'s,  after Andrei Zelevinsky brought to their attention the following two results.   Howlett proved that  any Coxeter group $\mathcal W$ is  infinite if and only if $c$ has infinite order ~\cite[Theorem 4.1]{h}.  Fomin and Zelevinsky proved the  following. Let $A$ be symmetrizable and bipartite, i.e., the set $\{1,\dots,n\}$ is a disjoint union of nonempty subsets $I,J$ and, for $h\ne l$,  $a_{hl}=0$ if either $h,l\in I$ or $h,l\in J$.  If $c=\prod_{i\in I}\sigma_i\prod_{j\in J}\sigma_j$, then ${\mathcal W}(A)$ is infinite if and only if the powers of $c$ are reduced words ~\cite[Corollary 9.6]{CAIV}.  The aforementioned result of ~\cite{kp} is a strengthening of the indicated results of Howlett and Fomin-Zelevinsky in the case  $\mathcal {W}=\mathcal {W}(A)$ where $A$ is symmetric.  A goal of  ~\cite{kp2} is to obtain the strengthening for any symmetrizable $A$, using properties of preprojective modules and of the set $\mathfrak S$.  Since there is a one-to-one correspondence between valued graphs and symmetrizable Cartan matrices ~\cite[p. 1]{dr}, one has to replace graphs with valued graphs and representations of quivers with representations of valued quivers.  Thus we lay a foundation for ~\cite{kp2}.

This paper continues  the study of combinatorial properties of   $\mathfrak S$  initiated in ~\cite{kt} and further developed in ~\cite{kp}.  These properties allow us to extend the main results of  ~\cite{kt} from representations of quivers to representations of valued quivers (Section \ref{S:preproj}), as well as to give new, more transparent proofs.  The rich combinatorics of $\mathfrak S$ is not fully understood and is useful for representation theory.  Our intention is to study it in the future. 

We now recall some facts, definitions, and notation, using
freely ~\cite{ars,bgp,dr}. A {\it graph} is a pair
$\Gamma=(\Gamma_{0},\Gamma_{1})$, where $\Gamma_{0}$ is the set of
vertices and $\Gamma_{1}$ is the set of edges, i.e.,  of two-element
subsets of $\Gamma_0$. Any subset $X\subset\Gamma_0$ determines a {\it full subgraph} of $\Gamma$ with the set of vertices $X$ and the set of edges consisting of all those two-element subsets $\{i,j\}\in\Gamma_1$ that satisfy $i,j\in X$. A {\it valuation} $\textrm{\textbf{b}}$ of a graph $\Gamma$ is a set of
nonnegative integers $\{b_{ij}\}$ for all pairs $i,j\in\Gamma_0$ where $b_{ii}=0$ and there exist nonzero natural numbers $d_i$ satisfying
\[d_ib_{ij}=d_jb_{ji},\ \ \textrm{for\ all}\ i, j\in\Gamma_0.\]
The pair $(\Gamma,\textrm{\textbf{b}})$ is a {\it valued graph}, and 
the above condition says that the matrix
$[b_{ij}]$ is symmetrizable. The valued graph $(\Gamma,\textrm{\textbf{b}})$ is {\it connected\,} if for all vertices $h\ne l$, there is a sequence $h,\dots,i,j,\dots,l$ in  $\Gamma_0$ such that $b_{ij}\ne0$ for each pair of subsequent vertices $i,j$.  Throughout the paper, $(\Gamma,\textrm{\textbf{b}})$ is a fixed finite
connected valued graph with $|\Gamma_{0}|>1$, where $|X|$ stands for the
cardinality of a set $X$.

An orientation, $\Lambda$, on $\Gamma$
consists of two functions $s:\Gamma_{1}\to\Gamma_{0}$ and
$e:\Gamma_{1}\to\Gamma_{0}$.  For an edge $a\in\Gamma_{1}$, $s(a)$
and $e(a)$ are the vertices incident with $a$, and they are called
the starting point and the endpoint of $a$, respectively; one writes
$a:s(a)\to e(a)$.  The ordered triple $(\Gamma,\textrm{\textbf{b}},\Lambda)$ is called a {\it valued quiver} and $a$
is then called an arrow of the quiver.  Any subset $X\subset\Gamma_0$ determines a {\it full subquiver} of $(\Gamma,\textrm{\textbf{b}},\Lambda)$ by taking the full subgraph of $\Gamma$ determined by $X$ and preserving  the valuation and orientation of each edge. Given a sequence of arrows
$a_{1}¥,\dots,a_{t}¥,\ t>0,$ satisfying $e(a_{i}¥)=s(a_{i+1}¥),\
0<i<t,$ one forms a path $p=a_{t}¥\dots a_{1}¥$ of length $t$ in
$(\Gamma,\textrm{\textbf{b}},\Lambda)$.  By definition, $s(p)=s(a_{1}¥),\ e(p)=e(a_{t}¥),$ so one
writes $p:s(p)\to e(p)$ and says that $p$ is a path from $s(p)$ to
$e(p)$.  By definition, for all $x\in\Gamma_{0}¥$ there is a unique
path of length $0$ from $x$ to $x$, denoted by $e_x$. A path $p$ of
length at least 1 is an oriented cycle if $s(p)=e(p)$. The set of
vertices of any valued quiver without oriented cycles (no finiteness
assumptions) acquires a structure of a partially ordered set (poset)
by putting $x\le y$ if there is a path from $x$ to $y$.  If $(\Gamma,\textrm{\textbf{b}},\Lambda)$
has no oriented cycles, we denote this poset by $(\Gamma_{0}¥,\La)$.  All orientations $\Lambda,\Theta$, etc., are
such that $(\Gamma,\textrm{\textbf{b}},\Lambda)$, $(\Gamma,\textrm{\textbf{b}},\Theta)$, etc., have no oriented cycles.

To define representations of a valued quiver $(\Gamma,\textrm{\textbf{b}},\Lambda)$, one has to choose a {\it modulation} $\mathfrak B$ of the valued
graph $(\Gamma,\textrm{\textbf{b}})$, which by definition is a set of division rings $\mathbf{k}_i$, $i\in\Gamma_0$,
together with a $\mathbf{k}_i-\mathbf{k}_j$-bimodule $_iB_j$ and a
$\mathbf{k}_j-\mathbf{k}_i$-bimodule $_jB_i$ for each edge
$\{i,j\}\in\Gamma_1$ such that

(i) there are $\mathbf{k}_j-\mathbf{k}_i$-bimodule isomorphisms
\[_jB_i\cong\textrm{Hom}_{\mathbf{k}_i}(_iB_j,\mathbf{k}_i)\cong\textrm{Hom}_{\mathbf{k}_j}(_iB_j,\mathbf{k}_j)\] and

(ii) dim$_{\mathbf{k}_i}({_iB_j})=b_{ij}$.

For the rest of the paper we denote by $\Gamma$ a valued graph with a fixed valuation $\textrm{\textbf{b}}$ and modulation $\mathfrak B$, denote by $\gl$ the  corresponding valued quiver with orientation $\La$, and assume that the division rings $\mathbf{k}_i$ are finite dimensional vector spaces over a common central subfield $k$ acting centrally on all bimodules $_iB_j$.  The latter assumption is sufficient for the applications that we have in mind.  However, the results of ~\cite{dr2}  imply that most of our considerations  hold without this assumption. Under the assumption, each   $_iB_j$  is a finite dimensional $k$-space, so setting $d_i=\dim_k\mathbf{k}_i$, we get $d_ib_{ij}=\dim_k({_iB_j})=\dim_k({_jB_i})=d_jb_{ji}$.

 A (left) representation $(V,f)$ of $\gl$ is a set of
finite dimensional left $\mathbf{k}_i$-spaces $V_i$, $i\in\Gamma_0$,
together with $\mathbf{k}_j$-linear maps
\[f_a:{_jB_i}\otimes_{\mathbf{k}_i}V_i\to V_j\] for all arrows $a:i\to
j$, and  morphisms of representations are defined in a natural way.  We obtain the category $\textrm{Rep}\gl$ of representations of the valued quiver $\gl$.

Putting
$\mathbf{k}=\prod_{i\in\Gamma_0}\mathbf{k}_i$ and viewing
$B=\underset{i\to j}\bigoplus{_jB_i}$ as a
$\mathbf{k}$-$\mathbf{k}$-bimodule where $\mathbf{k}$ acts on ${_jB_i}$ from the left via the projection $\mathbf{k}\to\mathbf{k}_j$ and from the right via the projection $\mathbf{k}\to\mathbf{k}_i$, one forms the tensor
ring $\textrm{T}(\mathbf{k},B)=\bigoplus_{n=0}^\infty B^{(n)}$ where
$B^{(n)}=B\otimes_\mathbf{k}\cdots\otimes_\mathbf{k} B$ is the
$n$-fold tensor product, and the multiplication is given by the
isomorphisms $B^{(n)}\otimes B^{(m)}\to B^{(n+m)}$ ~\cite[p.
386]{dr1}.  Since $\gl$ has no oriented cycles,
$\textrm{T}(\mathbf{k},B)$ is a finite dimensional $k$-algebra and we denote it by $k\gl$.  Let  $e_i\in\mathbf{k}$ be the $n$-tuple that has $1\in\mathbf{k}_i$ in the $i$th place and 0 elsewhere.  A
left $k\gl$-module $M$ is {\it finite dimensional} if
$\mathrm{dim}_{\mathbf{k}_i}e_iM<\infty$ for all $i$, which is equivalent to $\dim_k M<\infty$.
  We let f.d.$\,k\gl$
denote the category of finite dimensional left $k\gl$-modules.
The categories
Rep$\gl$ and f.d.$\,k\gl$ are equivalent ~\cite[Proposition 10.1]{dr1}
and we view the equivalence as an identification.  In this paper all $k\gl$-modules are finite dimensional.

Given a valued quiver $\gl$ and a vertex $x\in\Gamma_{0}$, let
$\sigma_{x}\Lambda$ be the orientation on $\Gamma$ obtained by
reversing the direction of each arrow incident with $x$ and
preserving the directions of the remaining arrows.  There results a
new valued quiver $\newspec$ (remember, the  valuation $\mathbf b$ and
modulation $\mathfrak B$ of the valued graph $\Gamma$ are fixed).  A vertex $x$ is a {\it sink}  if no arrow starts at
$x$.  For each sink $x$, the {\it reflection} functor
$F^{+}_{x}:\textrm{Rep}\gl\to\textrm{Rep}\newspec$ is defined
~\cite[pp. 15-16]{dr}, and we recall the definition for the
convenience of the reader.

Let $(V,f)\in\textrm{Rep}\gl$ and let $(W,g)=F^{+}_{x}(V,f)$.  Then
$W_y=V_y$ for all $y\neq x$, and $g_b=f_b$ for all those arrows
$b$ of $\newspec$ that do not start at $x$.  Let $a_i:y_i\to x,\
i=1,\dots,l$, be the arrows of $\gl$ ending at $x$. Then the
reversed arrows $a_i':x\to y_i,\ i=1,\dots,l$, are all the arrows of
$\newspec$ starting at $x$.  Consider the exact sequence
\[
0\to\Ker
h\overset{j}\to\overset{l}{\underset{i=1}\oplus}{_xB_{y_i}}\otimes_{\mathbf{k}_
{y_i}} V_{y_i}\overset{h}\to V_x
\]
of $\mathbf{k}_x$-spaces, where the map $h$ is induced by the maps
$f_{a_i}:{_xB_{y_i}}\otimes_{\mathbf{k}_ {y_i}}V_{y_i}\to V_x$. Then
$W_x=\Ker h$ and  each map
$g_{a_i'}:{_{y_i}B_x}\otimes_{\mathbf{k}_x}W_x\to W_{y_i}=V_{y_i}$ is
obtained from the map
$W_x\to{_xB_{y_i}}\otimes_{\mathbf{k}_{y_i}}W_{y_i}$ induced by $j$
using the following chain of isomorphisms of $k$-spaces
~\cite[pp. 14-15]{dr}.
\begin{eqnarray*}\textrm{Hom}_{\mathbf{k}_x}(W_x,{_xB_{y_i}}\otimes_{\mathbf{k}_{y_i}}W_{y_i})&\cong&
\textrm{Hom}_{\mathbf{k}_x}(W_x,\textrm{Hom}_{\mathbf{k}_{y_i}}({_{y_i}B_x},\mathbf{k}_{y_i})\otimes_{\mathbf{k}_{y_i}}W_{y_i})\\
&\cong&\textrm{Hom}_{\mathbf{k}_x}(W_x,\textrm{Hom}_{\mathbf{k}_{y_i}}({_{y_i}B_x},W_{y_i}))\\&\cong&
\textrm{Hom}_{\mathbf{k}_{y_i}}({_{y_i}B_x}\otimes_{\mathbf{k}_{x}}W_x,W_{y_i})\end{eqnarray*}
A sequence of vertices $S=x_{1},x_{2},\dots,x_{s},\ s\ge0,$ is called
{\it (+)-admissible} on $\gl$ if it either is  empty, or satisfies the
following conditions: $x_{1}$ is a sink with respect to $\Lambda$,
$x_{2}$ is a sink with respect to $\sigma_{x_{1}}\Lambda$, and so
on; sometimes we write $x_{1}x_{2}\dots x_{s}$ instead of
$x_{1},x_{2},\dots,x_{s}$. Recall that we denote by $\mathfrak S$ the set of
(+)-admissible sequences on $\gl$. If  $S=x_{1},\dots,x_{s}$ is in
$\mathfrak S$, we put
$\Lambda^{S}=\sigma_{x_{s}}\dots\sigma_{x_{1}}\Lambda$ and
$F(S)=F^{+}_{x_{s}}\dots F^{+}_{x_{1}}:\textrm{Rep}\gl\to
\textrm{Rep}(\Gamma, \Lambda^S)$. If the sequence $S$ consists of
distinct vertices and contains each vertex of the quiver, then
$F(S)=\Phi^{+}$ does not depend on the choice of $S$ and is called
the {\it Coxeter} functor ~\cite[p. 19]{dr}. For $S\in\mathfrak S$ we say
that $S$ \textit{annihilates} a module $M\in\mathrm{f.d.}\,k\gl$ if $F(S)(V,f)=0$,
where $(V,f)$ is the representation of $\gl$ identified with $M$.  In light
of this identification, we often write $F(S)M$ or $\Phi^+ M$.

A {\it source} is a
vertex of a quiver at which no arrow ends.  Replacing sinks with
sources, one gets similar definitions of a {\it reflection} functor  $F_{x}^{-}$, a
($-$){\it -admissible} sequence, and the {\it Coxeter} functor $\Phi^{-} $  ~\cite{dr}.

In ~\cite[p. 22]{dr}, the authors make the following definition.
\begin{defn}\label{intr1} A representation $(V,f)$ of $\gl$ is
{\it preprojective} if $(\Phi^{+})^{m}(V,f)=0$ for some integer
$m>0$.
\end{defn}
Definition \ref{intr1} is
equivalent to the following.
\begin{defn}\label{intr2} A module $M\in\mathrm{f.d.}\,k\gl$ is
{\it preprojective} if there exists an $S\in\mathfrak S$ that
annihilates it.
\end{defn}

We describe all $S\in\mathfrak S$
that annihilate a preprojective $k\gl$-module $M$ by proving that,  up
to a certain equivalence $\sim$, there exists a unique
{\it shortest} (+)-admissible sequence
$S_{M}$ that annihilates $M$ (Theorem \ref{shrtstsq}(a)), where an $S\in\mathfrak
S$ is a shortest sequence that annihilates $M$ if
$S$
annihilates $M$ but no proper subsequence of $S$ does.  Suppose now that $M$ is indecomposable. Then $S_{M}\in\mathfrak P$ (Theorem \ref{principal}) where
$\mathfrak P$ is the subset  of $\mathfrak
S$ consisting of the {\it principal} (+)-admissible sequences 
defined below in terms of the poset $(\Gamma_{0},\La)$ and
geometry of $\Gamma$, and $S_M$ determines $M$ uniquely up to isomorphism (Theorem \ref{shrtstsq}(d)).  If  $m$ is the smallest positive integer satisfying
$(\Phi^{+})^{m}M=0$, then $m=\nu+1$ where $\nu$ is a unique nonnegative integer for which
$(\Phi^{+})^{\nu}M=P$ is indecomposable projective; $P=P_{x}$ is determined up to
isomorphism by a unique $x\in\Gamma_{0}$; and $M\cong(\Phi^{-})^{\nu}P_{x}¥\cong(\TrD)^{\nu}P_{x}$. It easy to compute $S_{M}$
from $(\nu,x)=(\nu(M),x(M))$ and vice versa, and it is more efficient to compute $M$
from $S_{M}$ than from $(\nu,x)$ (Corollary \ref{invariants}).  If
$\gl$ is of infinite representation type, then
$\mathfrak P=\{S_{M}\,|\,M \mathrm{\ indecomposable\ preprojective}\}$ (Corollary \ref{posetisopreproj}(c)).   If $M_1,\dots,M_t$ are the nonisomorphic indecomposable summands of a preprojective module $M$, it is easy to compute $S_M$ in terms of $S_{M_1},\dots,S_{M_t}$ (Theorem \ref{shrtstsq}(c)).

The  preprojective
(connected) component of the Auslander-Reiten quiver of $k\gl$ is closely related to the translation quiver $\mathbb N\times(\Gamma,\La^{op})$, and if $\gl$ is of infinite representation type, the two coincide.  Recall that $\mathbb N\times(\Gamma,\La^{op})$, with $\mathbb
N$ being the set of nonnegative integers and $\La^{op}$ the opposite orientation of $\La$, is an infinite connected valued quiver that can be visualized as a disjoint union of countably many copies of the valued quiver $(\Gamma,\La^{op})$ where, for
each $i\in\mathbb N$, one draws
additional arrows 
starting at vertices of
$\{i\}\times(\Gamma,\La^{op})$ and ending at vertices of
$\{i+1\}\times(\Gamma,\La^{op})$; here  the valuation of new edges is assigned in a natural way and the translation is a left shift.  One of the reasons to study (+)-admissible sequences is that a significant part of the combinatorial structure of $\mathbb N\times(\Gamma,\La^{op})$ can be recovered from a simpler combinatorics of the set $\mathfrak S$, which has a natural poset
structure (up to the equivalence $\sim$): if $S,T\in\mathfrak S$, we set
$S\preccurlyeq T$ if $T\sim SS'$ where $S'$ is a (+)-admissible
sequence on $(\Gamma,\La^{S})$.  Since the translation quiver $\mathbb N(\Gamma,\La^{op})$ has no
oriented cycles, its set of vertices $\mathbb N\times\Gamma_{0}$
is a poset.  We prove that this poset is
isomorphic to $\mathfrak P$ viewed as a subposet of $\mathfrak S$ (Theorem \ref{posetiso}(a)).
A large class of valued quivers,
which is easy to describe combinatorially, is characterized by the fact that
$\gl$ with the valuation ignored coincides with the Hasse diagram of the poset
$(\Gamma_{0},\La)$. For these valued quivers,
the Hasse diagram of $\mathfrak P$ is the underlying quiver of the valued translation quiver
$\mathbb N(\Gamma,\La^{op})$ (Theorem \ref{posetiso}(b)), i.e., (+)-admissible sequences contain all
information about the preprojective component except for the valuation.

We now describe the content of the paper section by section. 
Section \ref{S:canform} presents the necessary definitions and results of ~\cite{kt, kp} concerning the combinatorics of the sets $\mathfrak S$ and $\mathfrak P$: the equivalence $\sim$, the partial order
$\preccurlyeq$, and a canonical
form and the lattice structure on the set $\mathfrak S$;  the
filters of the poset $(\Gamma_{0},\La)$ play a major role.   These considerations do not involve representation theory, valuation, or modulation of $\gl$, so most of the proofs are omitted.  Section 2 describes
 the properties of the shortest sequence $S_{M}$ associated to a preprojective module $M$, as
well as the connection between the preprojective component of the
Auslander-Reiten quiver of $\gl$ and the poset $\mathfrak P$.  

By duality, one can
study $(-)$-admissible sequences and the preinjective component of the valued
quiver, using ideals, instead of filters, of the poset
$(\Gamma_{0}¥,\La)$ and the same equivalence $\sim$.  We leave this to the reader.

The authors are grateful to Vlastimil Dlab and Claus Michael Ringel for helpful suggestions.

\section{Posets, Admissible Sequences, and Canonical Forms}\label{S:canform}


Throughout this section $\gl$ is a valued quiver without oriented cycles.
By definition, a (+)-admissible sequence on  $\gl$ depends
neither on the valuation \textbf{b} nor on the modulation $\mathfrak
B$.  Therefore the considerations
of ~\cite[Sections 1 and 2]{kt} and ~\cite[Section 2]{kp} apply and we quote, mostly without proofs, those results that are needed in the rest of the paper.  

We recall some notions about posets; see ~\cite{e}.  Let $(P,\leq)$
be a poset.  A subset $F\subset P$ is called a \textit{filter} if
whenever $x\in F$ and $y\geq x$, we have $y\in F$. We say that a
filter $F$ is generated by $X\subset P$ and write $F=\langle
X\rangle$ if $F=\{y\in P\,|\, y\geq x$ for some $x\in X\}$. If $F$ is
generated by a single element $x$, we call $F$ a \textit{principal
filter} and write $F=\langle x\rangle$.  For $x,y\in P$ we say that
$y$ covers $x$ and write $x\lessdot y$ if (i) $x<y$ and (ii)
$x<y'\le y$  implies $y'=y$. The Hasse diagram, ${\mathscr H}(P)$,
of $P$ is the quiver with the set of vertices $P$ and the set of
arrows that contains a single arrow $x\to y$ if and only if
$x\lessdot y$, and has no other arrows. 
For all
 $x,y\in\Gamma_{0}$, we set $x\leq y$ if there is a path from $x$
 to $y$ in $\gl$.  Since $\gl$ has no
 oriented cycles, 
 this turns $\Gamma_{0}$ into a poset, which we denote by
 $(\Gamma_{0},\Lambda)$.

\begin{defn}
If $S=x_1,\dots,x_s,\ s\ge0,$ is in $\mathfrak S$, we write $\La^S=\sigma_{x_s}\dots\sigma_{x_1}\La$ and, in particular,  $\La^{\emptyset}=\La$. The {\it support} of $S$,
$\mathrm{Supp}\, S$,  is the set of distinct vertices among $x_j$, $1\leq
j\leq s$.   As in ~\cite[Definition 2.1]{kp},  the {\it length} of $S$ is $\ell(S)=s$; the \textit{multiplicity} of $v\in\Gamma_0$ in $S$, $m_S(v)$,
is the (nonnegative) number of subscripts $j$ satisfying $x_j=v$. A sequence $K\in\mathfrak S$ is \textit{complete} if $m_K(v)=1$
for all $v\in\Gamma_0$. If $S=x_1,\dots,x_s$ and $T=y_1,\dots,y_t$
are $(+)$-admissible on $\gl$ and $(\Gamma,\Lambda^{S})$,
respectively, the concatenation of $S$ and $T$ is the sequence
$ST=x_1,\dots,x_s,y_1,\dots,y_t$.  If $K$ is complete,
$\Lambda^K=\Lambda$ so that if $m>0$, then $K^m$ denotes the
concatenation of $m$ copies of $K$ and $K^m\in\mathfrak S$.
\end{defn}

The following statement, which is ~\cite[Proposition 1.3]{kt}, relates
the elements of $\mathfrak S$ to filters of the poset $\poset$. In
particular, it tells us precisely when a subset of $\Gamma_0$ can be
realized as the support of a sequence
$S\in\mathfrak{S}$.

\begin{prop}\label{ondelta} Let
$\Omega\subset\Gamma_{0}$.  There exists a sequence
$S=x_1,\dots,x_s,\ s\ge0,$ in $\mathfrak S$ satisfying
$\mathrm{Supp}\,S=\Omega$ if and only if $\Omega$ is a filter of
$\poset$.  Moreover, if $\Omega\ne\emptyset$ is a filter, the sequence
$S=x_1,\dots,x_s$ can be chosen so that $x_1,\dots,x_s$ are
distinct.
\end{prop}

The following is ~\cite[Definition 1.2]{kt}.

\begin{defn}\label{defequiv} If a sequence
 $S=x_1,\dots,x_{i}¥,x_{i+1}¥,\dots,x_s,\ 0<i<s$, in $\mathfrak S$ has the property
  that no
 edge of $\Gamma$ connects $x_{i}¥$ with $x_{i+1}¥$, then
 $T=x_1,\dots,x_{i+1}¥,x_{i}¥,\dots,x_s$ is in $\mathfrak S$, and we set $SrT$.
 We denote by $\sim$ the equivalence relation that is a reflexive and
 transitive closure of the symmetric binary relation $r$.
\end{defn}

The above definition is motivated by the
fact that if distinct vertices $x$ and $y$ are both sinks in $\gl$, then $F^+_xF^+_y=F^+_yF^+_x$, as follows from the analog of  ~\cite[Lemma 1.2, proof of part 3)]{bgp} for representations of valued quivers.  Hence
$S\sim T$ implies $F(S)=F(T)$.

The following is ~\cite[Proposition 1.6]{kt}.

\begin{prop}\label{tfae} If $S,T\in \mathfrak S$ are nonempty and
consist of
distinct vertices, the following are equivalent.
\begin{itemize}
\item[(a)]$S\sim T$.
\item[(b)]$\mathrm{Supp}\,S=\mathrm{Supp}\,T$.
\item[(c)]$\Lambda^{S}=\Lambda^{T}$.
\end{itemize}
\end{prop}

The next result, which is ~\cite[Proposition 1.9]{kt}, produces a
canonical form in $\mathfrak S$.

\begin{prop}\label{canonicalform} Let $S\in \mathfrak S$ be nonempty.
\begin{itemize}
\item[(a)]We have $S\sim S_{1}S_{2}\dots S_{r}$ where, for
all $i$, $S_{i}$ consists of distinct vertices, and
$\mathrm{Supp}\,S_i=\mathrm{Supp}\,{S_{i}S_{i+1}\dots S_{r}}$. Further,
if $\mathrm{Supp}\,S_{i}\neq\Gamma_{0}$ then
$\mathrm{Supp}\,S_{i+1}\subsetneq\mathrm{Supp\,}S_{i}$.
\item[(b)]Let $T\sim T_{1}T_{2}\dots T_{q}\,$ be a nonempty sequence
in $\mathfrak S$ where, for all $j$, $T_{j}$ consists of distinct
vertices, and $\,\mathrm{Supp}\,T_{j}=\mathrm{Supp}\,T_{j}T_{j+1}\dots
T_{q}$. Then $S\sim T$ if and only if $r=q$ and $S_{i}\sim T_{i}$ on
$(\Gamma,\Lambda^{S_{1}\dots S_{i-1}}),\,i=1,\dots,r$.
\end{itemize}
\end{prop}

\begin{defn}\label{formsegmentsize} If
$S\sim S_{1}S_{2}\dots S_{r}$ is a nonempty sequence in $\mathfrak S$ and the $S_{i}$
satisfy the conditions of Proposition \ref{canonicalform}(a), we say
that $S_{1}S_{2}\dots S_{r}$ is the \textit{canonical form}, and $r$
is the \textit{size}, of $S$. If $S=S_{1}S_{2}\dots S_{r}$, we say
that $S$ is in the canonical form.  The size of the empty sequence is zero.
\end{defn}

By Proposition \ref{canonicalform}(b), the size of a nonempty
sequence is uniquely determined and each $S_{i}¥$ is unique up to
equivalence.

We quote ~\cite[Remark 2.1]{kp}.

\begin{rmk}\label{multcanon} In the setting of Proposition \ref{canonicalform}(a), if $v\in \Gamma_0$  then $v\in\mathrm{Supp}\,S_i$  if and  only if $m_S(v)\ge i$.
\end{rmk} 

According to ~\cite[Definition 1.5]{kt},  for each filter $F$ of $(\Gamma_{0},\Lambda)$, the \textit{hull} of
$F$ is the smallest filter of
$(\Gamma_{0}¥,\Lambda)$ containing $F$, as well as each vertex of
$\Gamma_{0}\setminus F$ that is connected by an edge to a vertex in
$F$.  The hull of $F$ is denoted by $H_{\Lambda}(F)$.

We quote ~\cite[Remark 2.2]{kp}.

\begin{rmk}\label{ConnectHull} If $F$ is a filter of  $(\Gamma_0,\La)$ and the full subgraph of $\Gamma$ determined   by $\mathrm{Supp}\,F$ is connected (for example,  if $F$ is a principal filter), then the full subgraph of $\Gamma$ determined   by $\mathrm{Supp}\,H_{\La}(F)$ is connected.
\end{rmk}

An effective way of constructing all possible (+)-admissible
sequences is given by the next statement, which is
~\cite[Proposition 1.11]{kt}.

\begin{prop}\label{tails} {\rm (a)} If $S=S_{1}S_{2}\dots S_{r}\in\mathfrak S$ is
a nonempty sequence in the canonical form  then, for all $i$,
$\mathrm{Supp}\,S_{i}$ is a filter of $(\Gamma_{0},\Lambda)$ and, for
$0<i<r$,
$H_{\Lambda}(\mathrm{Supp}\,S_{i+1})\subset\mathrm{Supp}\,S_{i}$.\vskip.04in

{\rm (b)} If $F_{1}\supset\dots\supset F_{r-1}\supset F_{r}$ is a
sequence of nonempty filters of $\,(\Gamma_{0},\Lambda)$ satisfying
$H_{\Lambda}(F_{i+1})\subset F_{i}$ for $\,0<i<r$, then there exists
a unique up to equivalence sequence $S_{1}S_{2}\dots
S_{r}\in\mathfrak S$ in the canonical form satisfying
$\mathrm{Supp}\,S_{i}=F_{i}$ for all $i$.
\end{prop}

We now introduce a partial order on the
set of equivalence classes of $\sim$, define
 the
subset $\mathfrak P$ of principal (+)-admissible sequences in $\mathfrak S$, and
relate the poset structure  of $\mathfrak P$ to the combinatorial structure of
the translation quiver ${\mathbb N}(\Gamma,\Lambda^{op})$.  We quote ~\cite[Definition 2.1]{kt}.

\begin{defn}If $S,T\in\mathfrak S$, we say that $S$ is
 a \textit{subsequence} of $T$ and write
$S\preccurlyeq T$ if $T\sim SU$ for some (+)-admissible sequence $U$.
\end{defn}

\begin{prop}\label{subsequence}\begin{itemize}\item[(a)] The relation $\preccurlyeq$ is a preorder.
\end{itemize}
Let $S,T\in\mathfrak S$.
\begin{itemize}
\item[(b)] We have $S\preccurlyeq T$ and $T\preccurlyeq S$ if and only if
$S\sim T$.
\item[(c)] If $S,T$ are nonempty and if $S_{1}\dots S_{r}$,
$T_{1}\dots T_{q}$ are their canonical forms, respectively,
then $S\preccurlyeq T$ if and only if $r\leq q$ and
$S_{i}\preccurlyeq T_{i}$ for $\,0<i\leq r$.
\item[(d)] $S\preccurlyeq T$ if and only if  for all $v\in\Gamma_0$, $m_S(v)\le m_T(v)$.
\end{itemize}
\end{prop}
\begin{proof} (a), (b), and (c) are proved in ~\cite[Proposition 2.1]{kt}.

(d) This is a direct consequence of (c), Proposition \ref{canonicalform}(a), and Remark \ref{multcanon}.
\end{proof}

We quote ~\cite[Corollary 2.4 and Proposition 2.5]{kp}.

\begin{cor}\label{multequiv} If ${S,T\in\mathfrak S}$, then ${S \sim T}$ if and only if for all $v\in\Gamma_0$, $m_S(v)= m_T(v)$.
\end{cor}

\begin{prop}\label{Prop:cancellation} Let ${S \in
\mathfrak{S}}$ and  let ${U, V}$ be (+)-admissible sequences on
$(\Gamma,\La^S)$.
\begin{itemize}

\item[(a)]${SU \preccurlyeq SV}$
if and only if ${U \preccurlyeq V}$.
\item[(b)]${SU \sim SV}$ if and only if
${U \sim V}$.

\end{itemize}
\end{prop}

\begin{proof}  Part (a) is an immediate consequence of Proposition \ref{subsequence}(d), and (b) follows directly from Corollary \ref{multequiv}.
\end{proof}

By Proposition \ref{subsequence}(b), the preorder $\preccurlyeq$
induces a partial order, which we denote by the same symbol, on the
set of equivalence classes of $\sim$ in $\mathfrak S$; when no
confusion arises we identify a sequence with its equivalence class. 
The poset  $\mathfrak S$ is a lattice ~\cite{kp}, as is demonstrated below.

The following is ~\cite[Definition 2.4]{kp}.
\begin{defn}\label{def:meetjoin}
Let ${S, T \in \mathfrak{S}}$ be nonempty and let
 ${S_{1}S_{2} \dots S_{r}},\ { T_{1}T_{2} \dots T_{q}}$ be their canonical forms, respectively, where without loss of generality we assume that ${r \leq q}$.  We set:
\begin{itemize}
\item[(a)] ${S \wedge T}$ to be the empty sequence if $\mathrm{Supp}\,{S} \cap \mathrm{Supp}\,{T}=\emptyset$; and if $\mathrm{Supp}\,{S} \cap \mathrm{Supp}\,{T}\ne\emptyset$, then  ${S \wedge T}$ is a (+)-admissible sequence with the
canonical form ${R_{1}R_{2} \dots R_{s}}$, where $s\le q,r$ is the largest integer satisfying ${\mathrm{Supp}\,{R_{i}} = \mathrm{Supp}\,{S_{i}} \cap \mathrm{Supp}\,{T_{i}}}\ne\emptyset$ for $0< i \leq s$. 

\item[(b)] ${S \vee T}$ to be a (+)-admissible sequence with the
canonical form ${R_{1}R_{2} \dots R_{q}}$, where $\mathrm{Supp}\,{R_{i}} = \mathrm{Supp}\,{S_{i}} \cup \mathrm{Supp}\,{T_{i}}$ for ${0< i \leq r}$, and  ${\mathrm{Supp}\,{R_{i}} = \mathrm{Supp}\,{T_{i}}}$ for
${r < i \leq q}$.
\end{itemize}

If $\emptyset$ is the empty sequence in $\mathfrak{S}$, then for all $S\in\mathfrak{S}$, we set  ${S \wedge \emptyset=\emptyset}$  and ${S \vee \emptyset=S}$.
\end{defn}

That $S \wedge T$ and $S
\vee T$ are in fact (+)-admissible sequences is contained in the proof of the
following statement, which is ~\cite[Proposition 2.6]{kp}.

\begin{prop}\label{lattice}
The poset of equivalence classes of $\sim$ in ${\mathfrak{S}}$ with the partial order
${\preccurlyeq}$ is a lattice where the operations of the greatest lower bound and the least upper bound  are ${\wedge}$ and ${\vee}$, respectively.
\end{prop}

\begin{proof} The intersection or union of two filters is always a filter.  If $F_1,F_2$ are filters of $(\Gamma_0,\La)$, then it is straight forward that $H_{\La}(F_1\cap F_2)\subset H_{\La}(F_1)\cap H_{\La}(F_2)$ and $H_{\La}(F_1\cup F_2)= H_{\La}(F_1)\cup H_{\La}(F_2)$.  Therefore, in view of Proposition \ref{tails}, we conclude that if $S,T\in\mathfrak S$, then 
${S \wedge T}$ and ${S \vee T}$ are in ${\mathfrak{S}}$.   It follows from  Proposition \ref{subsequence}, parts (c) and (d), that ${S \wedge T}$ and ${S \vee T}$  are the
greatest lower bound and the least upper bound, respectively, of $S$ and $T$.
\end{proof}

We quote  ~\cite[Theorem 2.7]{kp}.

\begin{thm}\label{Theorem:meetjoin}
Let ${S, T \in \mathfrak{S}}$.
\begin{itemize}
 \item[(a)]  ${S \sim \left(S
\wedge T \right)S'}, \ {T \sim \left(S \wedge T \right)T'}$ where $S',T'$ are (+)-admissible sequences on $(\Gamma,\Lambda^{S \wedge T})$ that are unique up to equivalence.
\item[(b)]${\mathrm{Supp}\,S' \cap \mathrm{Supp}\,T' = \emptyset}$.
\end{itemize}
\end{thm}

\begin{proof}  (a)  This is a direct consequence of Propositions \ref{lattice} and \ref{Prop:cancellation}(b).

(b)  By (a), we have $(S \wedge T) (S' \wedge T')\preccurlyeq S,T $, so Proposition \ref{lattice}  implies $(S \wedge T) (S' \wedge T')\preccurlyeq S\wedge T$ whence $S' \wedge T'=\emptyset$.  By Definition \ref{def:meetjoin}(a) and Proposition \ref{canonicalform}(a), ${\mathrm{Supp}\,S' \cap \mathrm{Supp}\,T' = \emptyset}$.
\end{proof}

The following is  ~\cite[Definition 2.2]{kt}.

\begin{defn}\label{tight} A sequence $S\in\mathfrak S$ is \textit{tight}
 if it is nonempty and its
canonical form $S_{1}S_{2}\dots S_{r}$ satisfies
$\mathrm{Supp}\,S_{i}=H_{\La}(\mathrm{Supp}\,S_{i+1})$ for $0<i<r$, and $S$ is \textit{principal} if it is tight and
$\mathrm{Supp}\,S_{r}$ is a principal filter. We denote by $\mathfrak
T$ ($\mathfrak P$) the set of tight (principal) sequences in
$\mathfrak S$; clearly, $\mathfrak P\subset\mathfrak
T\subset\mathfrak S$. By Proposition \ref{tfae}, a tight sequence is
uniquely determined by its size and the set $\mathrm{Supp}\,S_{r}$, so
we let $S_{r,x}$ denote the principal sequence of size $r$ with
$\mathrm{Supp}\,S_{r}=\langle x\rangle,\ x\in\Gamma_{0}$. Thus
$\mathfrak P=\{S_{r,x}\,|\,r\in\mathbb Z^{+},x\in\Gamma_{0}\}$ where
$\mathbb Z^{+}$ is the set of positive integers.
\end{defn}

We quote ~\cite[Remark 3.1]{kp}.

\begin{rmk}\label{ConnectPrinc}  It follows from Remark \ref{ConnectHull} that if $S\in\mathfrak P$, the full subgraph of $\Gamma$ determined by $\mathrm{Supp}\,S$ is connected.
\end{rmk}

The next statement, which is
~\cite[Corollary 2.3]{kt}, explains how to compare an arbitrary
sequence in $\mathfrak{S}$ with a tight one, and shows that the last
vertex of a principal sequence is uniquely determined.

\begin{prop}\label{tightend} Let $S\in\mathfrak S$ be
nonempty and $T\in\mathfrak T$ with canonical
forms $S_{1}\dots S_{r}$ and $T_{1}\dots T_{q}$, respectively.
\begin{itemize}
\item[(a)] We have $T\preccurlyeq S$ if and only if $q\le r$ and
$\mathrm{Supp}\,T_{q}\subset\mathrm{Supp}\,S_{q}$. If\,
$T=S_{q,x}\in\mathfrak P$ then $T\preccurlyeq S$ if and only if
$q\le r$ and $x\in\mathrm{Supp}\,S_{q}$.
\item[(b)] If $T=S_{q,x}=x_{1}¥,x_{2}¥,\dots,x_{t}¥$
then $x_{t}=x$.
\end{itemize}
\end{prop}

We now recall the notion of a translation quiver ~\cite[p. 47]{r}. 
If $\Delta=(\Delta_{0},\Delta_{1})$ is a locally finite graph with 
an orientation $\Theta$, the quiver $(\Delta,\Theta)$ is a {\it translation 
quiver} if it is equipped with a partially defined injective map 
$\tau:\Delta_{0}\to\Delta_{0}$, called the {\it translation} of 
$(\Delta,\Theta)$, 
such that for all $z\in\Delta_{0}$ in the domain of $\tau$ 
and all $y\in\Delta_{0}$ there is an arrow from $y$ to $z$ if and only if there is an arrow from $\tau z$ to $y$ (remember, in this paper no graph has multiple edges). In particular, the 
translation quiver ${\mathbb N}(\Gamma,\Lambda^{op})$ of the opposite quiver of 
$\gl$ is defined as follows.  The set of vertices of ${\mathbb 
N}(\Gamma,\Lambda^{op})$ is ${\mathbb N}\times\Gamma_{0}$, and each 
arrow $a:u\to v$ of $\gl$, which by definition is the only arrow  $u\to v$,  gives rise to two series of arrows, 
$(n,a^{\circ}):(n,v)\to(n,u)$ and 
$(n,a^{\circ})':(n,u)\to(n+1,v)$. The translation is defined 
by $\tau(n,u)=(n-1,u)$ for all $n>0$ and $u\in\Gamma_{0}$.  By 
construction, ${\mathbb 
N}(\Gamma,\Lambda^{op})$ is a locally finite quiver without oriented 
cycles, so ${\mathbb N}\times\Gamma_{0}$ is a poset. We note that since $\gl$ is a valued quiver, ${\mathbb N}(\Gamma,\Lambda^{op})$ is a {\it valued translation quiver} (see ~\cite[Sections VII.4 and VIII.1]{ars}).  However, we do not use the valuation on ${\mathbb N}(\Gamma,\Lambda^{op})$ because our method is to obtain information about the latter set using the combinatorics of $\mathfrak S$ and $\mathfrak P$, which are independent of the valuation or modulation on $\gl$.

We end this section by relating the Hasse
diagram of $\mathfrak P$ to ${\mathbb N}(\Gamma,\Lambda^{op})$.  Recall that if $a:x\to y$ is an arrow in a quiver, then a path $a_t\dots a_1:x\to y$ of length $t>1$ in the quiver is called a {\it bypass} of $a$. The following is
~\cite[Theorem 2.5]{kt}.

\begin{thm}\label{posetiso} Let $\mathfrak P$ be the set of
principal $(+)$-admissible sequences on $\gl$.
\begin{itemize}
\item[(a)] The map $\psi:\mathfrak P\to{\mathbb N}\times\Gamma_{0}$ given by
$\psi(S_{r,x})=(r-1,x)$
is an isomorphism of posets.
\item[(b)] Suppose no arrow in $\gl$ has a bypass.
Then $\psi$ induces an isomorphism of quivers
$\psi:{\mathscr H}(\mathfrak P)\to{\mathbb N}(\Gamma,\Lambda^{op}¥)$,
and
the map $S_{r,x}¥\mapsto S_{r-1,x}$,
$x\in\Gamma_{0}¥,\ r>1$, is a translation on ${\mathscr H}(\mathfrak P)$
that turns $\psi$ into an isomorphism of
translation quivers.
\end{itemize}
\end{thm}

\section{Preprojective Modules}\label{S:preproj}

Throughout this section $\gl$ is a quiver without oriented cycles with a fixed valuation $\mathbf b$ and modulation $\mathfrak B$.  We apply the combinatorial results of Section
\ref{S:canform} to the preprojective component of the Auslander-Reiten quiver.

\begin{defn}If $S=x_{1},\dots,x_{s},\ s>0,$ is in $\mathfrak S$, we
let $F(S)$ denote the composition of reflection functors
$F^{+}_{x_{s}}\dots F^{+}_{x_{1}}$; when $S=K$ is a complete
$(+)$-admissible sequence then $F(S)=\Phi^{+}$ is the (positive)
Coxeter functor ~\cite{dr}, and if $S=\emptyset$ then $F(S)$ is the identity functor on ${\mathrm {f.d.}}\,k\gl$.   We say that $S$ \textit{annihilates} a
$k\gl$-module $M$ if $F(S)M=0$; if, in addition, no proper subsequence
of $S$ annihilates $M$, then $S$ is a \textit{shortest} sequence
annihilating $M$.
\end{defn}

Recall that if $(V,f)\in\mathrm{Rep}\gl$, the {\it support} of $(V,f)$ is defined as $\mathrm{Supp}\,(V,f)=\{x\in\Gamma_0\,|\,V_x\ne0\}$.  If $M\in\mathrm{f.d.}\,k\gl$ and $(V,f)$ is the representation identified with $M$, then, by definition, $\mathrm{Supp}\,M=\mathrm{Supp}\,(V,f)$.
 
\begin{rmk}\label{support}   Let $M\in\mathrm{f.d.}\,k\gl$.   If $S\in\mathfrak S$ annihilates $M$, then $\mathrm{Supp}\,M\subset\mathrm{Supp}\,S_M$.  If $M$ is indecomposable, the full subgraph of $\Gamma$ determined by $\mathrm{Supp}\,M$ is connected.
\end{rmk}

For each $x\in\Gamma_0$, let $L_x\in\textrm{Rep}\gl$ be defined by $L_x=(V_i,f_a)$, where $V_i=0$ for $i\neq
x$, $V_x=\mathbf{k}_x$, and $f_a=0$ for all  arrows $a$.  That is, the
representations $L_x$ are the simple objects of $\textrm{Rep}\gl$. The
following is an analog of ~\cite[Corollary 1.1]{bgp}.

\begin{prop}\label{analog} Let $x_{1},\dots,x_{s},\ s>0,$ be a
$(+)$-admissible sequence on the valued quiver $\gl$.
\begin{itemize} \item[(a)] For any $i$ $(1\leq i\leq s)$,
$F_{x_1}^-\cdots F_{x_{i-1}}^-(L_{x_i})$ is either
$0$ or an indecomposable object in $\mathrm{Rep}\gl$.
\item[(b)] If $(V,f)\in\mathrm{Rep}\gl$ is indecomposable and
$F_{x_s}^+\cdots F_{x_{1}}^+(V,f)=0$, then for some $i$,
$(V,f)\cong F_{x_1}^-\cdots F_{x_{i-1}}^-(L_{x_i})$.
\end{itemize}
\end{prop}
\begin{proof}  The statement follows from ~\cite[Proposition 2.1]{dr} in the same way as ~\cite[Corollary 1.1]{bgp} follows from ~\cite[Theorem
1.1]{bgp}.
\end{proof}

The following result extends ~\cite[Theorem 3.1]{kt} and ~\cite[Theorem 3.4(a)]{kp} to representations of valued quivers.  For an integer $m>0$ and $N\in\mathrm{f.d.}\,k\gl$, we denote by $N^m$ the direct sum of $m$ copies of $N$.

\begin{thm}\label{shrtstsq} Let $M\in\mathrm{f.d.}\,k\gl$ be preprojective.
\begin{itemize} \item[(a)] There exists a unique up to
equivalence shortest sequence $S_M\in\mathfrak{S}$ annihilating
$M$.
\item[(b)]  If $M\cong N_1\oplus\dots\oplus N_s$ then each $N_i$ is preprojective  and $S_M=S_{N_1}\vee\dots\vee S_{N_s}$.  In particular, for all integers $m>0$, $S_{M^m}=S_M$.

\item[(c)]  If $M\cong M_1^{m_1}\oplus\dots\oplus M_t^{m_t}$ where the $M_i$'s are nonisomorphic indecomposable $k\gl$-modules and $m_i>0$ for all $i$, then $S_M=S_{M_1}\vee\dots\vee S_{M_t}$. 
\item[(d)]  If $M$ is indecomposable and $N\in\mathrm{f.d.}\,k\gl$ is 
indecomposable preprojective, then $S_N\sim S_M$ if and only if
$N\cong M$.
\end{itemize}
\end{thm}

\begin{proof} (a) Let $S,T$ be shortest sequences in $\mathfrak S$ annihilating $M$, where $\ell(S)\le\ell(T)$. To show that $S\sim T$, proceed by induction on $\ell(S)$. If $\ell(S)=0$, then $S=\emptyset$ whence $M=0$ and $T=\emptyset$.  Suppose now that $\ell(S)>0$ and that the statement holds for all preprojective $k(\Gamma,\Theta)$-modules $N$ (for all  orientations $\Theta$ without oriented cycles) and all pairs $S',T'$ of shortest (+)-admissible sequences on $(\Gamma,\Theta)$ annihilating $N$ where $\ell(S')\le\ell(T')$ and $\ell(S')<\ell(S)$.  Since $\ell(S)>0$, then $M\ne0$.

By Theorem \ref{Theorem:meetjoin}, $S\sim(S\wedge T)S'$ and $T\sim(S\wedge T)T'$ where $S',T'$ are (+)-admissible sequences  on $(\Gamma,\La^{S\wedge T})$ satisfying $\mathrm{Supp}\,S'\cap\,\mathrm{Supp}\,T'=\emptyset $.  It follows that $\ell(S')\le\ell(T')$. If $S\wedge T=\emptyset$ then $S\sim S', \, T\sim T',$ and $\mathrm{Supp}\,S\cap\mathrm{Supp}\,T=\emptyset $.  Let $(W,h)\in\mathrm{ Rep}\gl$ be identified with $M$.  Since $F(S)M=0$, if $W_i\ne0$ for some $i\in\Gamma_0$, then $i\in\mathrm{Supp}\,S$.  Since $\mathrm{Supp}\,S\cap\mathrm{Supp}\,T=\emptyset $, then $F(T)$ does not change any of the nonzero $\mathbf{k}_i$-spaces $W_i$, which exist because $M\ne0$.  We obtained a contradiction with $F(T)M=0$, so $S\wedge T\ne\emptyset$ whence $\ell(S')<\ell(S)$ and $S',T'$ are shortest (+)-admissible sequences on $(\Gamma,\La^{S\wedge T})$ annihilating the preprojective $k(\Gamma,\La^{S\wedge T})$-module $F(S\wedge T)M$.  By the induction hypothesis, we have $S'\sim T'$ whence $S\sim (S\wedge T)S'\sim(S\wedge T)T' \sim T$.

(b) Since every reflection functor is additive, each $N_i$ is  preprojective.  By (a), a sequence $S\in\mathfrak S$ annihilates $M$ if and only if $S_{N_i}\preccurlyeq S$ for all $i$.  Since $\mathfrak S$ is a lattice by Proposition \ref{lattice}, we have $S_M=S_{N_1}\vee\dots\vee S_{N_s}$. 

(c) This is an immediate consequence of (b).

(d)  This follows from Proposition \ref{analog}.
\end{proof}

As in ~\cite{kt}, we now show that if $M$ is an indecomposable  preprojective 
$k\gl$-module, then $S_{M}\in\mathfrak P$, and we begin with the case
when $M=P$ is projective.

\begin{lem}\label{strproj} Let $P_x\in\mathrm{f.d.}\,k\gl$ be the indecomposable projective
module associated with $x\in\Gamma_0$, let $(V,f)$ be the representation of $\gl$ identified with $P_x$, and let $x\ne z\in\Gamma_0$.
\begin{itemize}
\item[(a)] If  the set $\{a_i:y_i\to z,\ i=1,\dots,l\}$ of all arrows ending
at $z$ is not empty, then the map $h:\oplus_{i=1}^l
({_zB_{y_i}}\underset{\mathbf{k}_{y_i}}\otimes V_{y_i})\to V_z$
induced by the maps
$f_{a_i}:{_zB_{y_i}}\underset{\mathbf{k}_{y_i}}\otimes V_{y_i}\to
V_z$ is an isomorphism.
\item[(b)] If $z$ is a sink in $\gl$ and if $\,Q_x\in\mathrm{f.d.}\,k{(\Gamma,\sigma_z\La)}$
is the indecomposable projective module associated with
$x$, then $F_{z}^{+}(P_{x})\cong Q_{x}$.
\end{itemize}
\end{lem}

\begin{proof} (a) We recall the structure of $(V,f)$, see
~\cite[Section 10]{dr1} and ~\cite{dr2}.  For all $u\in\Gamma_0,$ denote by ${\mathcal W}^x_u$
the set of all paths  from $x$ to $u$ in $\gl$ and let $p\in{\mathcal W}^x_u$.   If
$p=b_t\cdots b_1, \ t>0,$ we set
$B_p={_{e(b_t)}B_{s(b_t)}}\underset{\mathbf{k}_{s(b_t)}}\otimes\cdots\underset{\mathbf{k}_{e(b_1)}}\otimes
{_{e(b_1)}B_{s(b_1)}}$, and if $u=x$ and $p=e_x$ is the trivial path, we set $B_p=\mathbf{k}_x$.    Then
$V_z=\underset{p\in{\mathcal W}^x_z}\oplus B_p$; note that   $V_z=0$  if
${\mathcal W}^x_z=\emptyset$.  To describe the map
associated to an arrow $y\to z$, say, to $a_1:y_1\to z$, we note
first that
$$V_z=(\underset{p=a_1q}\oplus
B_p)\oplus(\underset{p\ne a_1q}\oplus
B_p)=(\underset{q\in{\mathcal W}^x_{y_1}}\oplus
({_zB_{y_1}}\underset{\mathbf{k}_{y_1}}\otimes
B_q))\oplus(\underset{p\ne a_1q}\oplus B_p),$$ while
$V_{y_1}=\underset{q\in{\mathcal W}^x_{y_1}}\oplus B_q$.  The
function $f_{a_1}:  {_zB_{y_1}}\underset{\mathbf{k}_{y_1}}\otimes
V_{y_1}\to V_z$ maps its domain onto the first summand of its
codomain via the usual isomorphism
${_zB_{y_1}}\underset{\mathbf{k}_{y_1}}\otimes({\underset{q\in{\mathcal
W}^x_{y_1}}\oplus} B_q)\to \underset{q\in{\mathcal
W}^x_{y_1}}\oplus ({_zB_{y_1}}\underset{\mathbf{k}_{y_1}}\otimes
B_q)$.  It is now clear that $h$ is an isomorphism.

(b) Let $(U,j)$ and $(W,g)$ be the representations of
$(\Gamma,\sigma_z\La)$ identified with $Q_x$ and $F_z^+(P_x)$,
respectively.  Let $z\ne y\in\Gamma_0$. Since $z$ is a sink in $\gl$ and a source in $(\Gamma,\sigma_z\La)$, 
a path from $x$ to $y$ in $(\Gamma,\sigma_z\Lambda)$  is a path from $x$ to $y$ in $\gl$, and vice versa. Thus, $U_y=V_y$ for all $y\neq z$,
and $j_a=f_a$  for all arrows $a$ not ending at
$z$.  Since $F_z^+$ affects only the  space at
 $z$ and the maps into this space, $U_y=W_y$  for
all $y\neq z$, and $j_a=g_a$ for all arrows $a$ not ending at $z$.  Since $z$
is a source in $(\Gamma,\sigma_z\Lambda)$, there is no path from $x$
to $z$ in that quiver, whence $U_z=0$.   It remains to show that
$W_z=0$.  Since  $z$ is a sink in $\gl$, $x\ne z$, and the graph
$\Gamma$ is connected, the set of arrows stopping at $z$ is not
empty.  Hence the map $h$ of part (a)  is an isomorphism, and
$W_z=\Ker h=0$.
\end{proof}

\begin{prop}\label{projectivesequence}If $P_{x}$ is the indecomposable
projective $k\gl$-module associated with $x\in\Gamma_{0}$, then
$S_{P_{x}}\sim S_{1,x}$.
\end{prop}

\begin{proof} By Propositions \ref{ondelta} and \ref{tfae}, there exists
a unique  up to equivalence (+)-admissible sequence
$S=x_{1},\dots,x_{s}$ that consists of distinct vertices and
satisfies $\{x_{1},\dots,x_{s}\}=\langle x\rangle$.  We first show
by induction on $s$ that $S$ annihilates $P_{x}$.  When $s=1$ this
follows from Proposition \ref{analog}. Suppose $s>1$ and the
statement holds for all orientations $\Theta$ on $\Gamma$ without oriented
cycles and all indecomposable projective $k(\Gamma,\Theta)$-modules associated with vertices $w$ satisfying $|\langle w\rangle|<s$. Since $s>1$, then
$x<x_{1}$ in $(\Gamma_{0},\Lambda)$, and in
$(\Gamma_{0},\sigma_{x_{1}}\Lambda)$ we have $\langle
x\rangle=\{x_{2},\dots,x_{s}\}$. By the induction hypothesis, the
sequence $x_{2},\dots,x_{s}$ annihilates $Q_{x}$, the indecomposable
projective $k(\Gamma,\sigma_{x_{1}}\Lambda)$-module associated with
$x$. Since $x<x_1$ and $x_1$ is a sink in $\gl$, Lemma \ref{strproj}(b)
says that $F_{x_{1}}^{+}(P_{x})\cong Q_{x}$, so $S$ annihilates
$P_{x}$. To show that no proper subsequence of $S$ annihilates
$P_{x}$, let $(V,f)\in\mathrm{Rep}\gl$ be identified with $P_{x}$
 and note that if $y\ge x$ in $(\Gamma_{0},\Lambda)$, then
$V_y\neq 0$. Since $V_y$ may be changed by 
$F_{z}^{+}$ only if $z=y$, any sequence annihilating $P_{x}$ must contain $y$.
\end{proof}

We need the following combinatorial statement whose necessity is ~\cite[Proposition 3.6]{kp}.

 \begin{prop}\label{princseq} Let $S=x_1,\dots,x_s,\ s>1,$ be in $\mathfrak S$ and set  $T=x_2,\dots,x_s$.  Then $S\in\mathfrak P$  if and only if  the full subgraph of $\Gamma$ determined by $\mathrm{Supp}\,S$ is connected and $T$ is a principal (+)-admissible sequence on $(\Gamma,\sigma_{x_1}\La)$.  
 \end{prop} 
 \begin{proof} We only have to prove the sufficiency.  Parts of the proof are similar to the proofs of ~\cite[Proposition 3.6 and Theorem 4.5]{kp}.
 
For a given graph $\Gamma$ and orientation $\La$, the sets $\mathfrak{P}$ and ${\mathfrak S}$ depend neither on the valuation $\mathbf b$ nor on the modulation $\mathfrak B$.  Hence, without loss of generality, we may assume for the rest of this proof that $\gl$ is an ordinary, not valued, quiver in which at least one of the arrows has multiplicity greater than $1$: since $\Gamma$ is a connected graph with more than one vertex, $\gl$ has at least one arrow.  Then the finite dimensional path algebras $k\gl$ and $k(\Gamma,\sigma_{x_1}\La)$ are of infinite representation type (see ~\cite{bgp}), and the results of ~\cite{kt} apply.
 
 Since $k(\Gamma,\sigma_{x_1}\La)$ is of infinite representation type, ~\cite[Corollary 3.8, parts (a) and (c)]{kt} says that $T\sim S_N$ for some indecomposable preprojective $k(\Gamma,\sigma_{x_1}\La)$-module $N$, and   ~\cite[VIII Proposition 1.14]{ars} says that $N$ is not a preinjective module, hence, not a simple injective  module.  By ~\cite[Theorem 1.1, part (2)]{bgp}, $M=F_{x_1}^-N$ is an indecomposable $k\gl$-module and $N\cong F_{x_1}^+M$, so that $S$ annihilates $M$. Therefore $M$ is preprojective, $S_M\preccurlyeq S$, and  ~\cite[Theorem 3.5]{kt} says that $S_M\in\mathfrak P$.  To show that $S\in\mathfrak P$, it suffices to prove that $x_1\in\mathrm{Supp}\,S_M$.  For if the latter is true, then $S_M\sim y_1,\dots,y_t$ where $y_1=x_1$ and 
 \[F_{y_t}^+\dots F_{y_2}^+(N)\cong F_{y_t}^+\dots F_{y_2}^+F_{y_1}^+(F_{x_1}^-N)=F(S_M)M=0, \]
 whence $y_2,\dots,y_t$ is a (+)-admissible sequence on $(\Gamma,\sigma_{x_1}\La)$ that annihilates $N$, so that $\ell(S_M)-1\ge\ell(T)=\ell(S)-1$ and $\ell(S_M)\ge\ell(S)$.  Since $S_M\preccurlyeq S$, then $S_M\sim S$.
 
 If $x_1\not\in\mathrm{Supp}\, S_M$ then $x_1\in\mathrm{Supp}\, U$, where $S\sim S_MU$, and $x_1$ is a sink in $(\Gamma,\La^{S_M})$ because $\mathrm{Supp}\, S_M$, being a filter of $(\Gamma_0,\La)$, contains  no $v\in\Gamma_0$ satisfying $v\le x_1$.  By ~\cite[Lemma 1.7]{kt}, for all $v\in\mathrm{Supp}\, S_M$, no arrow connects $v$ and $x_1$ whence $S_Mx_1\sim x_1S_M$ on $\gl$.  Therefore $S_M$ is a (+)-admissible sequence on $(\Gamma,\sigma_{x_1}\La)$ and we have $0=F^+_{x_1}(F(S_M)M)=F(S_M)(F^+_{x_1}M)\cong F(S_M)N$.  Hence $S_N\preccurlyeq S_M$ so that $s-1\le\ell(S_M)$, which implies $s-1=\ell(S_M)$ and $S\sim S_Mx_1\sim x_1S_M$.  Then the full subgraph of $\Gamma$ determined by $\mathrm{Supp}\, S$ is disconnected, a contradiction. 
\end{proof}

Although  the statement of Proposition \ref{princseq} does not involve representation theory, our proof uses representations of quivers.  We know a  purely combinatorial proof, but it is much longer and more technical than the  one given above. 

The next result is an extension of ~\cite[Theorem 3.5]{kt} to representations of valued quivers.

\begin{thm}\label{principal} If $M\in\mathrm{f.d.}\,k\gl$ is indecomposable preprojective,
$S_{M}$ is a principal $(+)$-admissible sequence.
\end{thm}

\begin{proof} Since $M\ne0$, then $S_M=x_1,\dots,x_s, \ s>0,$ and we proceed by induction on $s$.  The case $s=1$ is trivial, so let $s>1$ and suppose that the theorem holds for all orientations $\Theta$ on $\Gamma$ without oriented cylces and all indecomposable preprojective  $k(\Gamma,\Theta)$-modules $N$ satisfying $\ell(S_N)<s$.  Since $s>1$, $N=F^+_{x_1}M$  is an indecomposable preprojective  $k(\Gamma,\sigma_{x_1}\La)$-module and $S_N=x_2,\dots,x_s$.  By the induction hypothesis, $S_N$ is a principal (+)-admissible sequence on $(\Gamma,\sigma_{x_1}\La)$.  In view of Proposition \ref{princseq}, to prove that $S_{M}\in\mathfrak P$, it suffices to show that the full subgraph of $\Gamma$ determined by $\mathrm{Supp}\,S_M$ is connected.  

Assume, to the contrary, that the subgraph is disconnected.  Since  $S_N$ is a principal (+)-admissible sequence, Remark \ref{ConnectPrinc} says that the full subgraph of $\Gamma$ determined by $\mathrm{Supp}\,S_N$ is connected, whence $\mathrm{Supp}\,S_M=\mathrm{Supp}\,S_N\cup\{x_1\}$  where $x_1\not\in\mathrm{Supp}\,S_N$ and, moreover, no edge of $\Gamma$ connects $x_1$ to a vertex in $\mathrm{Supp}\,S_N$.  It follows that $S_M=x_1S_N\sim S_Nx_1$ so that $S_N\in\mathfrak S$.  According to Remark \ref{support}, the full subgraph of $\Gamma$ determined by $\mathrm{Supp}\,M$ is connected and $\mathrm{Supp}\,M\subset\mathrm{Supp}\,S_M$.  Then either $\mathrm{Supp}\,M=\{x_1\}$ or $\mathrm{Supp}\,M\subset\mathrm{Supp}\,S_N$.  In the former case, $M\cong  L_{x_1}$ whence $S_M=x_1$, which contradicts $s>1$. In the latter case, $0=F(S_M)M=F^+_{x_1}(F(S_N)M)$ implies $F(S_N)M=0$ because $x_1\not\in\mathrm{Supp}\,S_N$, which contradicts that $S_M$ is the shortest sequence annihilating $M$.
\end{proof}

\begin{cor}\label{invariants} Let $M\in\mathrm{f.d.}\,k\gl$ be indecomposable and
satisfy $(\Phi^{+})^{\nu}M\cong P_{x}$, where $\nu\in{\mathbb N}$
and $P_{x}$ is the indecomposable projective $k\gl$-module associated
with $x\in\Gamma_{0}$.
\begin{itemize}
\item[(a)] $S_{M}\sim S_{\nu+1,x}$.
\item[(b)] If $S_{M}=x_{1},\dots,x_{s}$ then $x_{s}=x$ and $M\cong
F_{x_{1}}^{-}\dots F_{x_{s-1}}^{-}(L_{x})$ where $L_{x}$ is the
simple projective
$k{(\Gamma,\sigma_{x_{s-1}}\dots\sigma_{x_{1}}\La)}$-module
associated with $x$.
\end{itemize}
\end{cor}

\begin{proof} (a) By Theorem \ref{principal}, $S_{M}\sim S_{r,y}$. We
have $S_{M}\preccurlyeq K^{r-1}S_{1,y}\preccurlyeq K^{r}$ by
Propositions \ref{tails} and \ref{subsequence}, so
$F(S_{1,y})((\Phi^{+})^{r-1}M)=(\Phi^{+})^{r}M=0$. Since
$(\Phi^{+})^{r-1}M\ne0$ by Theorem \ref{shrtstsq}(a), then
$(\Phi^{+})^{r-1}M\cong P_{x}$ and $\nu=r-1$ (see
\cite[Proposition 2.4(i)]{dr}). Since $S_{1,y}$ annihilates $P_{x}$
then $S_{1,x}\preccurlyeq S_{1,y}$ by Proposition
\ref{projectivesequence}. Since $K^{r-1}S_{1,x}$ annihilates $M$
then $S_{M}\preccurlyeq K^{r-1}S_{1,x}$, whence $S_{1,y}\preccurlyeq
S_{1,x}$ and $S_{1,y}\sim S_{1,x}$ in light of Proposition
\ref{subsequence}. Using Proposition \ref{tfae}, we get $x=y$.

(b) This is an easy consequence of (a), Corollary \ref{tightend}(b),
and Proposition \ref{analog}.
\end{proof}

 In order to apply our results to the
 preprojective component  of $\gl$, we recall
 some definitions and facts from ~\cite {ars,r}.
  If $X\in\mathrm{f.d.}\,k\gl$ is
 indecomposable, let
 $[X]$ be the isomorphism class of $X$. If $Y\in\mathrm{f.d.}\,k\gl$ is
 indecomposable, a path of length $m>0$ from $X$ to $Y$ is a
 sequence of
 nonzero nonisomorphisms $X=A_{0}\to\dots\to A_{m}=Y,$ where
 $A_{i}\in\mathrm{f.d.}\,k\gl$ is
 indecomposable for all $i$. By definition, there
  exists a path of
 length zero from $X$ to $X$.  One writes $[X]\prec [Y]$ if there
  exists a path of positive length from $X$ to $Y$.

  The preprojective component of $\gl$,
 $\tilde{\mathscr P}\gl$, is a locally finite connected valued translation
 quiver whose set of vertices, $\tilde{\mathscr P}\gl_{0}$, consists
 of the isomorphism classes of
 indecomposable
 preprojective $k\gl$-modules.  If $X,Y\in\mathrm{f.d.}\,k\gl$ are
 indecomposable, there is an arrow $[X]\to[Y]$ 
 if and only if  there exists an irreducible map $X\to Y$ (remember, we disregard the valuations of arrows).  The translation is defined by
 $[X]\mapsto[\DTr X]=[\Phi^{+}¥X]$ for all nonprojective $X$.  If
 $X,Y$ are indecomposable, $Y$ is preprojective, and $X=A_{0}¥\to\dots\to
 A_{m}¥=Y,\ m>0,$ is a path from $X$ to $Y$, then $[X]\ne[Y]$ and
 $A_{i}¥$ is preprojective for all $i$.  It follows that the reflexive closure
 $\preccurlyeq$
 of the transitive binary relation $\prec$ is a partial order on
 $\tilde{\mathscr P}\gl_{0}¥$.
 Moreover, $[X]\prec [Y]$ if and only if there is a finite
 sequence of
 irreducible morphisms $X=B_{0}¥\to\dots\to B_{n}¥=Y$, where $n>0$
 and $B_{j}¥$ is
 indecomposable preprojective for all $j$.
 
 We finish the paper by extending ~\cite[Proposition 3.7 and Corollary 3.8]{kt} to representations of valued quivers.  Consider the map
 $\phi:\tilde{\mathscr P}\gl\to{\mathbb N}(\Gamma,\Lambda^{op}¥)$
 defined on the vertices by
$\phi([L])=(\nu,x)=(\nu(L),x(L))$, where $x$ is the vertex of $\gl$
 associated with the indecomposable projective module
 $(\Phi^{+}¥)^{\nu}¥L$, and defined on the arrows in a natural way
 ~\cite[VIII Proposition 1.15]{ars}.

\begin{prop}\label{components}
 \begin{itemize}
 \item[(a)] The map
 $\phi:\tilde{\mathscr P}\gl\to{\mathbb N}(\Gamma,\Lambda^{op}¥)$ is
  a full embedding
 of translation quivers whose restriction
 $\phi:\tilde{\mathscr P}\gl_{0}¥\to{\mathbb N}\times\Gamma_{0}¥$ is
 an injective morphism of posets.
 \item[(b)] The map $\phi$ is an isomorphism when $\gl$ is of
 infinite representation type.
 \item[(c)] The image of $\phi$ is an
 ideal of $\,{\mathbb N}\times\Gamma_{0}¥$, i.e., if
 $\,[M]\in\tilde{\mathscr P}\gl_{0}¥$ and
 $(l,u)\le\phi([M])$, then there
 exists
 an indecomposable preprojective  $k\gl$-module
 $L$ with $\phi([L])=(l,u)$.
 \item[(d)] Given an $[M]\in\tilde{\mathscr P}\gl_{0}¥$, the map $\phi$
 induces a bijection between the set of paths in $\tilde{\mathscr P}\gl$
 ending at $[M]$ and the set of paths in ${\mathbb N}(\Gamma,\Lambda^{op}¥)$
 ending at $\phi([M])$.
 \end{itemize}
 \end{prop}

 \begin{proof} (a) and (b) These are ~\cite[VIII Propositions 1.15 and
 1.16]{ars}.

 (c) This is an easy consequence of the following obvious statement.  If
 $0\to A\to B\to C\to0$ is an almost split
 sequence of finitely generated modules over a hereditary artin algebra where $B$
  has a nonzero
 injective direct summand, then $C$ is injective.

 (d) This is an immediate consequence of (a) and (c).
  \end{proof}

 We now obtain a  module-theoretic version of Theorem \ref{posetiso}.

 \begin{cor}\label{posetisopreproj}
\begin{itemize}
\item[(a)] The map $\chi:\tilde{\mathscr P}\gl_{0}¥\to\mathfrak P$ given by
$[L]\mapsto S_{L}$ is an injective morphism of posets.
\item[(b)] If each arrow $x\to y$ is the only path from $x$ to $y$ in $\gl$,
then
the map $\chi$ induces a full embedding
$\chi:\tilde{\mathscr P}\gl\to{\mathscr H}(\mathfrak P)$ of translation
quivers,
where $S_{r,x}¥\mapsto S_{r-1,x}¥$,
$x\in\Gamma_{0}¥,\ r>1$, is the translation on ${\mathscr H}(\mathfrak P)$.
\item[(c)]  If
$\gl$ is of infinite representation type, the map $\chi$ in (a) and
in (b)
is an
isomorphism.
\end{itemize}
 \end{cor}
 \begin{proof}  This is an immediate consequence of Theorems
 \ref{shrtstsq},
 \ref{principal}, and
 \ref{posetiso}, together with Proposition \ref{components}.
 \end{proof}

\end{document}